\documentclass[journal]{IEEEtran}

\ifCLASSINFOpdf
\usepackage[pdftex]{graphicx}
 \DeclareGraphicsExtensions{.pdf,.jpeg,.png}
\else
 \usepackage[dvips]{graphicx}
 \DeclareGraphicsExtensions{.eps}
\fi
\usepackage[cmex10]{amsmath}
\usepackage{array}
\usepackage{mdwmath}
\usepackage{mdwtab}
\usepackage{amsfonts}
\usepackage{amsthm}
\usepackage{tikz}
\theoremstyle{plain}
   \newtheorem{thm}{Theorem}
   
   \newtheorem{lemma}[thm]{Lemma}
   
   \theoremstyle{definition}
   
   \newtheorem{defn}[thm]{Definition}
   \newtheorem{example}[thm]{Example}
   \theoremstyle{remark}
   
   \newtheorem{remark}[thm]{Remark}

   \newcommand{\refthm}[1]{Theorem~\ref{#1}}
   
   \newcommand{\reflemma}[1]{Lemma~\ref{#1}}
   
   \newcommand{\refdefn}[1]{Definition~\ref{#1}}
   \newcommand{\refremark}[1]{Remark~\ref{#1}}
   
\newcommand{\Hom}{\operatorname{Hom}}
   \newcommand{\Spec}{\operatorname{Spec}}
\newcommand{\orb}{\operatorname{orb}}
\newcommand{\C}{\operatorname{H}}
\newcommand{\Res}{\operatorname{Res}}
\newcommand{\F}{\mathbb F_q^*}
\usepackage{fixme}
\usepackage{amscd}
\usepackage{amsmath,amssymb}
\usepackage[all,line]{xy}
\begin{document}

\title{Quantum Codes from Toric Surfaces}

\author{Johan~P.~Hansen}% <-this % stops a space
\thanks{Johan P. Hansen is with the Department
of Mathematics, Aarhus University, Aarhus, Denmark, e-mail: matjph@imf.au.dk.}% <-this % stops a space
\thanks{Part of this work was done while visiting Institut de Math\'ematiques de Luminy, MARSEILLE, France.
I thank for the hospitality shown to me.}% <-this % stops a space

\markboth{Preprint March~2012}%
{Shell \MakeLowercase{\textit{et al.}}: Bare Demo of IEEEtran.cls for Journals}

\maketitle

\begin{abstract}

A theory for constructing quantum error correcting codes from 
Toric surfaces by the Calderbank-Shor-Steane method is presented. In particular we study the method on toric Hirzebruch surfaces.

The results are obtained by constructing a dualizing differential form for the toric surface and by using the cohomology and the intersection theory of toric varieties.
\end{abstract}


\section{Introduction}
In  \cite{MR1749454} and \cite{MR1953195}  the author developed methods to construct linear error correcting codes from toric varieties and derive the code parameters using the cohomology and the intersection theory on toric varieties. This method is generalized in section \ref{toric} to construct linear codes suitable for constructing quantum codes by the Calderbank-Shor-Steane method. Essential for the theory is the existence and the application of a dualizing differential form on the toric surface.

A.R. Calderbank \cite{Calderbank19961098}, 
P.W. Shor \cite{MR1450603} and 
A.M. Steane \cite{Steane19992492} produced stabilizer codes from linear codes containing their dual codes.

These two constructions are merged to obtain results for toric surfaces in section \ref{final}. Similar merging has been done for algebraic curves  with different methods by A. Ashikhmin, S. Litsyn and M.A. Tsfasman in \cite{tsfasman}.

\hfill Johan P. Hansen
\hfill March 20, 2012

\subsection{Notation}
\begin{itemize}
\item $\mathbb F_q$ -- the finite field with $q$ elements of characteristic $p$.
\item $\F$ -- the invertible elements in $\mathbb F_q$.
\item $k=\overline{\mathbb F_q}$ -- an algebraic closure of $\mathbb F_q$.
\item $M \simeq \mathbb Z^2$ a free $\mathbb Z$-module of rank 2.
\item $\square \subseteq M_{\mathbb R}= M \otimes_{\mathbb Z}{\mathbb R}$ -- an integral convex polytope.
\item $X=X_{\square}$ -- the toric surface associated to the polytope $\square$.
\item $T=T_N =U_0 \subseteq X$ -- the torus.
\item $S = |D_1| \cap |D_2| \subseteq X(\mathbb F_q)$ -- the intersection of the supports of the divisors $D_1$ and $D_2$.
\item $\omega_X$ the sheaf of differential forms on $X$.
\end{itemize}

\section{The method of toric varieties}\label{toric}
For the general theory of toric varieties we refer to \cite{MR2810322}, \cite{MR1234037} and  \cite{MR922894}.
Here we will be using toric surfaces and we recollect some of their theory.

\subsection{Toric surfaces and their cohomology}

Let $M$ be an integer lattice $M \simeq \mathbb Z^2$. Let
$N=\Hom_\mathbb Z(M,\mathbb Z)$ be the dual lattice with canonical $\mathbb Z$ -
bilinear pairing $ <\quad,\quad>: M \times N \rightarrow \mathbb Z.$ Let
$M_{\mathbb R}= M \otimes_{\mathbb Z}{\mathbb R}$ and $N_{\mathbb
R}= N\otimes_{\mathbb Z}{\mathbb R}$ with canonical  $\mathbb R$ -
bilinear pairing $ <\quad,\quad>: M_{\mathbb R} \times N_{\mathbb
R} \rightarrow  \mathbb R. $

Given a 2-dimensional integral convex polytope $\square$ in $M_{\mathbb R}$. The support function
$
h_{\square}: N_{\mathbb R} \rightarrow \mathbb R
$
is defined as
$
h_{\square}(n):= {\rm inf}\{<m,n> | \, m \in \square\}
$
and the polytope $\square$ can be reconstructed from the support function
\begin{equation}\label{support}
\square_{h} = \{m \in M |\, <m,n> \,\geq \,h(n) \quad \forall n \in N \}.
\end{equation}

The support function $h_{\square}$ is piecewise linear in the sense
that  $N_{\mathbb R}$ is the union of
a non-empty finite collection of strongly convex polyhedral cones in $N_{\mathbb R}$ such that
$h_{\square}$ is linear on each cone. A fan is a collection $\Delta$ of
strongly convex polyhedral cones in $N_{\mathbb R}$
such that every face of $\sigma \in \Delta $ is contained in $\Delta $ and $\sigma \cap \sigma'  \in \Delta$ for all
$\sigma , \sigma' \in \Delta $.

The {\it normal fan} $\Delta$ is the coarsest  fan such that $h_{\square}$ is linear
on each $\sigma \in \Delta$, i.e. for all $\sigma \in \Delta$ there exists  $l_{\sigma} \in M$ such that
\begin{equation}
h_{\square}(n) = <l_{\sigma},n> \quad \forall n \in \sigma.
\label{linear}
\end{equation}

The 1-dimensional cones $\rho \in \Delta$ are generated by unique primitive elements
$n(\rho) \in N \cap \rho$ such that
$\rho =\mathbb R_{\geq 0} n(\rho)$.

Upon refinement of the normal fan, we can assume that two
successive pairs of $n(\rho)$'s generate the lattice and we obtain
{\it the refined normal fan}, which will be the fan we will be using for the the rest of the present paper.

The 2-dimensional {\it algebraic torus} $T_N \simeq 
k^* \times k^* $ is defined by $T_N:= \Hom_\mathbb Z(M,k^*)$. The multiplicative
character $\mathbf e(m),\, m \in M$ is the homomorphism $\mathbf e(m): T \rightarrow k^*$ defined
by $\mathbf e(m)(t) = t(m)$ for $t \in T_N$. Specifically, if $\{n_1,n_2\}$ and $\{m_1,m_2\}$ are dual $\mathbb Z$-bases of $N$ and $M$ and we denote
$u_j:= \mathbf e(m_j),\,j=1,2$, then we have an isomorphism $T_N \simeq k^* \times k^* $ sending $t$ to $(u_1(t),u_2(t))$.
For $m=\lambda_1 m_1 +\lambda_2 m_2$ we have
\begin{equation}\label{e}
\mathbf e(m)(t)=u_1(t)^{\lambda_1}u_2(t)^{\lambda_2}.
\end{equation}
The {\it toric surface} $X_{\square}$ associated to the refined
normal fan $\Delta$ of $\square$ is
\begin{equation*}
X_{\square} = \cup_{\sigma \in \Delta} U_{\sigma}
\end{equation*}
where $U_{\sigma}$ is the $k$-valued points of the affine scheme  $\Spec(k[\mathcal S_{\sigma}])$, i.e.,
morphisms $u : {\mathcal S}_{\sigma} \rightarrow k$ with $u(0)=1$ and $
u(m+m')= u(m)u(m')\ \forall m,m' \in \mathcal S_{\sigma}$,
where $\mathcal S_{\sigma}$ is the additive subsemigroup of $M$
\begin{equation*}
\mathcal S_{\sigma}=\{m \in M | <m,y> \geq 0 \ \forall y \in \sigma\}.
\end{equation*}

The {\it toric surface} $X_{\square}$ is irreducible, non-singular
and complete under the assumption that we are working with the refined normal fan.
 If $\sigma, \tau \in \Delta$ and $\tau$ is a face of
$\sigma$, then $U_{\tau}$ is an open subset of $U_{\sigma}$.
Obviously $\mathcal S_{0}=M$ and $U_{0}=T_N$ such that the algebraic
torus $T_N$ is an open subset of $X_{\square}$.

$T_N$ {\it acts algebraically} on $X_{\square}$. On $u \in U_{\sigma}$ the action of $t \in T_N$ is obtained as
\begin{equation*}
(tu)(m):=t(m)u(m) \quad  \mathrm{for\ }m \in \mathcal S_{\sigma}\ ,
\end{equation*}
such that $tu \in U_{\sigma}$ and $U_{\sigma}$ is $T_N$-stable.
The orbits of this action are in one-to-one correspondance with
$\Delta$. For each $\sigma \in \Delta$, let
$$\orb(\sigma):=\{u:M\cap \sigma \rightarrow k^* | u \text{
is a group homomorphism}\}.$$ Then $\orb(\sigma)$ is a $T_N$ orbit
in $X_{\square}$. Define $V(\sigma)$ to be the closure of
$\orb(\sigma)$ in $X_{\square}$.

A $\Delta$-linear support function $h$ gives rise to a polytope $\square$ as above and an associated Cartier
divisor 
\begin{equation*}
D_h=D_{\square}:= -\sum_{\rho \in \Delta (1)} h(n(\rho))\,V(\rho)\ ,
\end{equation*}
where $\Delta (1)$ is the 1-dimensional cones in $\Delta$.
In particular
\begin{equation*}
D_m={\rm div}(\mathbf e(-m)) \quad m \in M.
\end{equation*}

\begin{lemma} \label{cohomology} Let  $h$ be a $\Delta$-linear support function
with associated convex polytope $\square$ and Cartier divisor $D_h=D_{\square}$. The vector space
${\rm H}^0(X,\it O_X(D_h))$ of global sections of $O_X(D_{\square})$, i.e., rational functions $f$ on
$X_{\square}$ such that
${\rm div}(f) + D_{\square} \geq 0$ has dimension $\#(M \cap \square)$ and has
$
\{\mathbf e(m) | m \in M \cap \square\}
$
as a basis.
\end{lemma}
\subsection{Intersection theory on a toric surface}
\label{mindist} 

For a $\Delta$-linear support function $h$ and a 1-dimensional
cone $\rho \in \Delta (1)$ we will determine the intersection
number $(D_h;V(\rho))$ between the Cartier divisor $D_h$ and
$V(\rho)) =\mathbb P^1$. This number is obtained in \cite[Lemma
2.11]{MR922894}. The cone $\rho$ is the common face of two 2-dimensional
cones $\sigma', \sigma'' \in \Delta (2)$. Choose primitive
elements $n', n'' \in N$ such that
\begin{align*}
n'+n''& \in \mathbb R \rho\\
\sigma' + \mathbb R \rho &= \mathbb R_{\geq 0} n' + \mathbb R \rho\\
\sigma'' + \mathbb R \rho &= \mathbb R_{\geq 0} n'' + \mathbb R \rho
\end{align*}
\begin{lemma}
\label{inter}For any $l_{\rho} \in M$, such that $h$ coincides
with $l_{\rho}$ on $\rho$, let $\overline{h} = h-l_{\rho}$. Then
\begin{equation*}
(D_h;V(\rho))= -(\overline{h}(n')+\overline{h}(n'').
\end{equation*}
\end{lemma}
In the 2-dimensional non-singular case let $n(\rho)$ be a
primitive generator for the 1-dimensional cone $\rho$. There
exists an integer $a$ such that
\begin{equation*}
n'+n''+a n(\rho)=0,
\end{equation*}
$V(\rho)$ is itself a Cartier divisor and the above gives the
self-intersection number
\begin{equation*}
(V(\rho);V(\rho))=a.
\end{equation*}
\fxnote{the 1-dim cones generate ...}
More generally the self-intersection number of a Cartier divisor
$D_h$ is obtained in \cite[Prop. 2.10]{MR922894}.
\begin{lemma} \label{self}Let  $D_h$   be a  Cartier divisor  and let  $\square_h$ be the polytope associated to
$h$. Then
\begin{equation*}
(D_h;D_h)= 2 \, {\rm vol}_{2}(\square_h),
\end{equation*}
where ${\rm vol}_{2}$ is the normalized Lesbesgue-measure.
\end{lemma}
\subsection{The support of the codes}
The toric codes are obtained from evaluating certain rational functions in a suitable set $S$ of $\mathbb F_q$-rational points on toric varieties, being the intersection of two ample divisors on $X$. 
\begin{defn}\label{support}
For $i=1,2$ let $I_i, J_i \subseteq \F$ with $I_1\cap J_2=I_2\cap J_1 = \emptyset$ and introduce the two rational functions
\begin{equation*}\label{funk}
F_i=\prod_{\psi \in I_1}(e(m_1)-\psi)^{n_{1,\psi}}\ \prod_{\psi \in J_1}(e(m_2)-\psi)^{n_{2,\psi}}\ ,
\end{equation*}
where the integer exponents satisfy $n_{1,\psi} \geq 1$ and $n_{2,\psi}\geq 1$.

For $i=1,2$, let $D_i= (F_i)_0$ be their divisor of zeroes, $\vert D_i \vert$ be their support and $U_i=X\backslash \vert D_i \vert$ their complement. It is important to note that the supports and their complement are independent of the choice of the exponents $n_{1,\psi} \geq 1$ and $n_{2,\psi}\geq 1$.

Finally let the support set of the code be $S =\vert D_i \vert \cap \vert D_2 \vert= U_1 \cup U_2 \subseteq \F \times \F$.
\end{defn}
\begin{remark}\label{dynamic}
As a set $S=I_1 \times J_2 \cup I_2 \times J_1 \subseteq \F \times \F$ with $\# S = \# I_1 \cdot \# J_2 + \# I_2 \cdot \# J_1$  elements, but it is important to  have in mind, that $S \subseteq \F \times \F$ is realized as the support of the intersection of two divisors in many different ways, namely one for each choice of the exponents $n_{1,\psi} \geq 1$ and $n_{2,\psi}\geq 1$.
\end{remark}

\subsection{Toric evaluation codes}\label{dimension}
We start by exhibiting the toric codes as evaluation codes supported on $S$.
\begin{defn}\label{toriccode}
For each $t \in T \simeq k^{*} \times k^{*}$, we
evaluate the rational functions in ${\rm H}^0(X,\it O_X(D_{\square}))$
\begin{eqnarray*}
    {\rm H}^0(X,\it O_X(D_{\square}))& \rightarrow & k\\
    f & \mapsto & f(t).
\end{eqnarray*}
Let ${\rm H}^0(X,\it O_X(D_{\square}))^{{\rm Frob}}$ denote the rational
functions in ${\rm H}^0(X,\it O_X(D_{\square}))$ that are invariant under
the action of the Frobenius, that is functions that are $\mathbb F_q$-linear
combinations of the functions $\mathbf e(m)$ in
(\ref{e}).

Evaluating in all points in $S$, we obtain the code
$C_{S,\square} \subset (\mathbb F_q)^{\# S}$ as the image
\begin{eqnarray*}
    {\rm H}^0(X,\it O_X(D_h))^{{\rm Frob}}& \rightarrow & C_{S,\square} \subset (\mathbb F_q)^{\# S}
 \\
    f & \mapsto & (f(t))_{t \in T(\mathbb F_q)}
\end{eqnarray*}
and the generators of the code is obtained as the image of the basis
\begin{equation*}
\mathbf e(m) \mapsto (\mathbf e(m)(t))_{t \in S}.
\end{equation*}
as in (\ref{e}).
\end{defn}
To estimate the parameters we have two bound the number of points in the support $S$, where the rational functions in ${\rm H}^0(X,\it O_X(D_{\square}))^{{\rm Frob}}$ evaluates to zero. 

The support $S$ is stratified by the intersections with the zeros of $\mathbf e(m_1)-\psi$, where $\psi \in I_1 \cup I_2 \subseteq \F$. A rational function $f$ can either vanish identically on a stratum or have a finite number of zeroes along the stratum.

\subsubsection{Identically vanishing} Assume that $f$ is identically zero
along precisely $a$ of these strata. As $\mathbf
e(m_{1})-\psi$ and $\mathbf
e(m_{1})$ have the same divisors of poles, they have equivalent
divisors of zeroes, so
\begin{equation*}
(\mathbf
e(m_{1})-\psi)_{0} \sim (\mathbf
e(m_{1}))_{0}.
\end{equation*}
Therefore
\begin{equation*}
{\rm div}(f) + D_{\square} -a (\mathbf
e(m_{1}))_{0}\geq 0
\end{equation*}
or equivalently
\begin{equation*}
f \in {\rm H}^0(X,\it
O_X(D_{\square}-a (\mathbf
e(m_{1}))_{0}).
\end{equation*}
Depending on the polytope $\square$ this gives an upper bound for the number $a$, using \reflemma{cohomology}.

\subsubsection{Vanishing in a finite number of points}
On any of the $\# I_1 \cup \# I_2-a$ other strata
the number of zeroes of $f$ is according to \cite{MR1866342} at most
the intersection number
\begin{equation}\label{intersection}
(D_{\square}-a (\mathbf e(m_{1}))_{0};(\mathbf
e(m_{1}))_{0}).
\end{equation}
 This number can be calculated using
\reflemma{inter} and  \reflemma{self}.

The above gives a method to construct toric codes from surfaces and obtain their precise parameters, this was done by the author in four cases in \cite{MR1953195}. 

\begin{example} (Hirzebruch surfaces).
Let $d, e, r$ be  positive integers and let $\square$ be the polytope
in $M_{\mathbb R}$ with vertices $(0,0),(d,0),(d,e+rd),(0,e),$ see
Figure \ref{hirzpolytop} and with (refined) normal fan as in Figure \ref{fan2}.
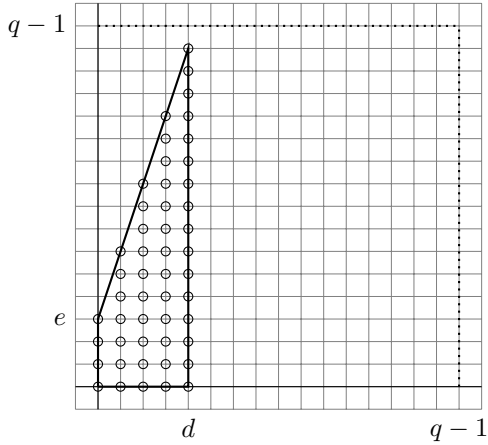
\begin{figure}
\begin{center}
\begin{tikzpicture}[scale=0.3]
\draw[step=1cm,gray,very thin] (-1,-1) grid (17,17);
\draw (-1,0)--(17,0);
\draw (0,-1)--(0,17);
\draw[thick,dotted](0,16)--(16,16);
\draw[thick,dotted](16,0)--(16,16);
\draw[thick] (0,0)--(4,0)--(4,15)--(0,3)--(0,0);
\draw (-1,3) node[anchor=east] {$e$};
\draw (-1,16) node[anchor=east] {$q-1$};
\draw (4,-1) node[anchor=north] {$d$};
\draw (16,-1) node[anchor=north] {$q-1$};
\foreach \x in {0,1,2,3} \draw (0,\x) circle (0.2cm);
\foreach \x in {0,1,2,3,4,5,6} \draw (1,\x) circle (0.2cm);
\foreach \x in {0,1,2,3,4,5,6,7,8,9} \draw (2,\x) circle (0.2cm);
\foreach \x in {0,1,2,3,4,5,6,7,8,9,10,11,12} \draw (3,\x) circle (0.2cm);
\foreach \x in {0,1,2,3,4,5,6,7,8,9,10,11,12,13,14,15} \draw (4,\x) circle (0.2cm);
\end{tikzpicture}
\end{center}
\caption{The polytope of \refthm{saetning4} is the polytope with
vertices $(0,0),(d,0),(d,e+rd),(0,e).$} \label{hirzpolytop}
\end{figure}
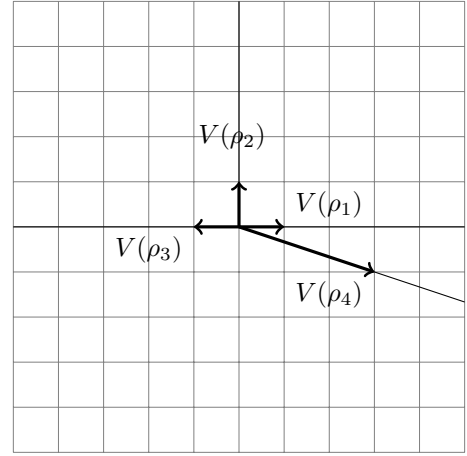
\begin{figure}
\begin{center}
\begin{tikzpicture}[scale=0.6]
\draw[step=1cm,gray,very thin] (-5,-5) grid (5,5);
\draw (-5,0)--(5,0);
\draw (0,0)--(0,5);
\draw (0,0)--(5,-5/3);
\draw[->, very thick] (0,0)--(0,1);
\draw[->, very thick] (0,0)--(-1,0);
\draw[->, very thick] (0,0)--(1,0);
\draw[->, very thick] (0,0)--(3,-1);
\draw (0.8,2) node[anchor=east] {$V(\rho_2)$};
\draw (-2,0) node[anchor=north] {$V(\rho_3)$};
\draw (2,1) node[anchor=north] {$V(\rho_1)$};
\draw (2,-1) node[anchor=north] {$V(\rho_4)$};
\end{tikzpicture}
\end{center}
\caption{The normal fan and the 1-dimensional cones of the polytope in Fig.
\ref{hirzpolytop} that gives rise to a Hirsebruch surface.} \label{fan2}
\end{figure}

\begin{table}\label{snittal}
\begin{center}
\begin{tabular}{c|cccc}
   & $V(\rho_1)$ & $V(\rho_2)$ & $V(\rho_3)$ & $V(\rho_4)$ \\ \hline
  $V(\rho_1)$& $-r$ & 1 & 0 & 1 \\
  $V(\rho_2)$ & 1 & 0 & 1 & 0 \\
  $V(\rho_3)$ & 0 & 1 & r & 1 \\
  $V(\rho_4)$ & 1 & 0 & 1 & 0 \\
\end{tabular}
\caption{The intersection numbers for the four 1-dimensional cones of the fan of the Hirzebruch surface}
\end{center}
\end{table}

From the Hirzebruch surfaces with $I_1=J_2=\F \times \F$ and $I_2=J_1=\emptyset$, we obtain using the above method  the following theorem.
\begin{thm} \label{saetning4} 
 Assume that  $d<q-1$, that  $e<q-1$ and
that $e+rd<q-1.$ The toric code $C_{\F\times \F,\square}$ has length equal to
$(q-1)^2,$ dimension equal to $\# (M \cap
\square)=(d+1)(e+1)+r\frac{d(d+1)}{2}$ (the number of lattice
points in $\square$) and the minimal distance is equal to
$
\mathrm{Min} \{(q-1-d)(q-1-e), (q-1)(q-1-e-rd)\}.
$
\end{thm}
 D. Joyner \cite{MR2142431}  has done extensive calculations on among others these toric codes. R. Joshua and R. Akhtbar \cite{MR2747728}
have obtained results on a different kind of toric codes that appear to be related to the dual of the present codes.
\end{example}

\section{Codes from toric surfaces containing their dual code}\label{merging}

\subsection{Differential forms and residues}
The residue theorem is obtained in \cite{MR0435074} over $\mathbb C$, however the theorem and its various forms are essential and for
completeness we present general proofs here. Throughout $\Res_P(\omega)$ means the local Grothendieck residue see, e.g., \cite{MR2464546} and \cite{MR759943}. For residues on toric varieties we also refer to \cite{MR1458757}.

\begin{thm}[Residue theorem - general form]\label{general}
Let $X$ be a complete smooth algebraic surface and let $\omega_X$ be the sheaf of differential 2-forms on $X$. Let $U_1, U_2$ be two open subsets of $X$ such that $X\backslash (U_1 \cup U_2) = S$ is a finite set of points. Then
\begin{itemize}
\item[i)]
Let $\omega \in \omega_X(U_1\cup U_2)=\C^0(U_1\cup U_2, \omega_X)$ be any 2-form on $X$ with no poles on $U_1\cup U_2$, then
$
\sum_{P \in S} \Res_P(\omega)=0\ .
$
\item[ii)] For any $(w_{\bullet}) \in \bigoplus_{P \in S} k$ with $\sum_{P \in S} w_P=0$, there exists an 
 $\omega \in \omega_X(U_1\cup U_2)=\C^0(U_1\cup U_2, \omega_X)$, such that $\Res_P(\omega)=w_P$ for all $P \in S$.
\end{itemize}
\end{thm}
\begin{IEEEproof}
The  \v{C}ech resolution 
${\omega_X}_{\vert{U_1}}\coprod {\omega_X}_{\vert{U_2}}\rightarrow {\omega_X}_{\vert{U_1\cup U_2}}$
of the sheaf ${\omega_X}_{\vert{U_1\cup U_2}}$, obtained from the two open sets $U_1$ and $U_2$, gives that a 2-form $\omega$ on $X$ without poles on $U_1\cup U_2$ defines a class $[\omega] \in  \C^1(U_1\cup U_2, \omega_X)$ and that every class has such a representation.

As $\C^2(X,\omega_X)\cong k$ and $\C^i(U_1\cup U_2,\omega_X)=0$ for $i \geq 2$ by Serre duality, relative cohomology gives the exact sequence
\begin{equation*}
\xymatrix@C=6pt@R=7pt{
\C^1(U_1\cup U_2,\omega_X) \ar[r]\ar@2{-}[dd]&\C_S^2(X,\omega_X)\ar[r]\ar[d]^{\cong} &\C^2(X,\omega_X)\ar[dd]^{\cong}\\
& \bigoplus_{P \in S}\C_P^2(X,\omega_X)\ar[d]^{\cong}& \\
\C^1(U_1 \cup U_2,\omega_X) \ar[r]_{[\omega] \mapsto {\oplus_{P \in S}}\Res_p(\omega)}^{\Res}
& \bigoplus_{P \in S}k\ar[r]_{\sum} & k
}
\end{equation*}
Then $\sum_{P \in S} \Res_P(\omega)= \Res([\omega])$ and the claims follows from exactness of the last sequence.
\end{IEEEproof}
In the above form there is no restrictions on the polar behavior as long as there are no poles on $U_1\cup U_2$, however it is possible to prove the theorem in a stronger form.

For a divisor $D$ on $X$ the sheaf of differential forms $\omega_X(D)$, is the sheaf with $\omega_X(D)(U)=\{\eta \in \omega(U) \vert (\eta)+D \geq 0\mathrm{\ on\ } U\}$ on open sets $U \subseteq X$. Its global sections  $\C^0(X,\omega(D))$ are the differential forms $\omega$ with $(\omega)+D\geq 0$. 

\begin{thm}[Residue theorem - special form]\label{special}
Let $X$ be a complete smooth algebraic surface and let $\omega_X$ be the sheaf of differential 2-forms on $X$. 
For $i=1,2$, let $D_i$  be ample and effective divisors on $X$ with support $\vert D_i\vert$ and with complement
 $U_i=X\backslash\vert D_i \vert$. Assume that $X\backslash (U_1 \cup U_2) = \vert D_1 \vert \cap \vert D_2 \vert= S$ is a finite set of points.
\begin{itemize}
\item[i)] For any $\omega \in \C^0(X,\omega(D_1+D_2))$, we have that
$\sum_{P \in S} \Res_P(\omega)=0\ .$
\item[ii)] For any $(w_{\bullet}) \in \bigoplus_{P \in S} k$ with $\sum_{P \in S} w_P=0$, there exists an 
$\omega \in \C^0(X,\omega(D_1+D_2))$, such that $\Res_P(\omega)=w_P$ for all $P \in S$.
\end{itemize}
\end{thm}
\begin{IEEEproof}
The morphism of sheaves
\begin{equation*}
\omega_X(D_1) \oplus \omega_X(D_2)\rightarrow \omega_X(D_1+D_2)\
\end{equation*}
is injective with cokernel $ j_{*}(\omega_{\vert U_1 \cap U_2})$, where $j$ is the open immersion of $U_1 \cup U_2 \hookrightarrow X$.
The associated long exact cohomology sequence gives a surjection
$
\C^0(X,\omega(D_1+D_2)) \rightarrow \C^1(U_1\cup U_2,\omega_X)
$ as $\C^1(X,\omega(D_1)) = \C^1(X,\omega(D_2)) = 0$ by the assumption on ampleness of the divisors.

The proof now follows as in the above proof of Theorem \ref{general}.
\end{IEEEproof}
\fxnote{bemærkninger om ample}
\subsection{Dualizing differential form of a toric code}
We want to exhibit a differential form $\omega_0$ on $X$ with poles restriced to the points in the support $S =\vert D_i \vert \cap \vert D_2 \vert= U_1 \cup U_2\subseteq \F \times \F$,  where $D_i= (F_i)_0$ are divisors of zeroes of the functions defined in \refdefn{support} and $\vert D_i \vert$ are their support and $U_i=X\backslash \vert D_i \vert$ their complement. Besides that we want the differential form $\omega$ to vanish at the divisor $2 D_{\square}$.
\begin{defn} A differential form $\omega_0 \in \C^0(X,\omega(D_1+D_2-2 D_{\square}))$ is called a dualizing form for the toric code and we will call the set 
\begin{equation*}
R =\{P \in S \vert \Res_P(\omega_0) \neq 0\} \subseteq S
\end{equation*}
its restricted support.
\end{defn}

This existence of a dualizing form for the toric code is obtained in two steps utilizing the representations of the set $S$ as the intersection of the supports of various ample divisors.

\begin{thm}Assume that the support of the toric code is the intersection of the support of two ample divisors as i Definition \ref{support}.
Assume that we can choose large exponents $n_{1,\psi} \geq 1$ and $n_{2,\psi}\geq 1$, such that $L(D_1+D_2-2 D_{\square})\neq \emptyset$.
Then there exists a dualizing form for the toric code of \refdefn{toriccode}  with ample divisors $D_1$ and $D_2$.
\end{thm}
\begin{IEEEproof}
For $i=1,2$, let $D_i= (F_i)_0$ be their divisor of zeroes, $\vert D_i \vert$ be their support and $U_i=X\backslash \vert D_i \vert$ their complement. Assuming that we can choose the exponents $n_{1,\psi} \geq 1$ and $n_{2,\psi}\geq 1$ such that $D_i$ are ample \refthm{special}
gives that for any $(w_{\bullet}) \in \bigoplus_{P \in S} k$ with $\sum_{P \in S} w_P=0$, there exists an 
$\omega \in \C^0(X,\omega(D_1+D_2))$, such that $\Res_P(\omega)=w_P$ for all $P \in S$.

In order to find a differential form vanishing at the divisor $2 D_{\square}$ we note that the support of the divisors $D_1$ and $D_2$ and their complement is independent of the choice of the exponents $n_{1,\psi} \geq 1$ and $n_{2,\psi}\geq 1$, see \refremark{dynamic}.  An $\omega$ constructed as above is in the corresponding $\C^0(X,\omega(D_1+D_2))$ for larger values of the exponents and the corresponding divisors $D_1$ and $D_2$ are still ample.

Choose large exponents $n_{1,\psi} \geq 1$ and $n_{2,\psi}\geq 1$, such that $L(D_1+D_2-2 D_{\square})\neq \emptyset$ and let $F \neq 0$ in $ L(D_1+D_2-2 D_{\square})$ using \reflemma{cohomology}.

The corresponding divisors $D_1$ and $D_2$ and the differential form $\omega_0=F \omega \in \C(X,\omega(D_1+D_2-2 D_{\square}))$ are the desired entities.
\end{IEEEproof}
\subsection{Toric codes contained in their dual codes}\label{final}
For a linear code $C \subseteq \mathbb F_q^n$ and $w \in (\F)^n$, we let the $w$-dual be the code
\begin{equation*}
C_w^{\perp}=\{x \in \mathbb F_q^n \vert \sum_1^n w_i x_i y_i =0 \ \forall y \in C\} \subseteq \mathbb F_q^n\ .
\end{equation*}
If $w=(1,1,\dots,1)$ then $C_w^{\perp}$ is the usual dual code $C^{\perp}$.
In the notation above, we have the following theorem.
\begin{thm}\label{dual} Let $C_\square$ be a toric code. Assume that the support $S$ of the toric code is the intersection of the support of two ample divisors as i Definition \ref{support}. Let $\omega_0 \in \C^0(X,\omega(D_1+D_2-2 D_{\square}))$ 
be a dualizing differential for the code with ample divisors $D_1$ and $D_2$ and let $R \subseteq S$ be the corresponding resticted support of the code. The evaluation code $C_{R,\square}$ obtained from evaluating functions in $L(D_{\square})$ at points in the restricted support $R$ satisfies $C_{R,\square} \subseteq  (C_{R,\square})_{w_{\bullet}}^{\perp}$, where $w_{P}= \Res(\omega_0)_{P}$ for $P \in R$.
\end{thm}
\begin{IEEEproof}For any $f,g \in L(D_{\square})$ and $\omega_0 \in \C(X,\omega(D_1+D_2-2 D_{\square}))$, we have that 
$f g \omega_0 \in \C(X,\omega(D_1+D_2)$. By \refthm{special} 
$\sum_{P \in R} f(P) g(P) \Res(\omega_0)_P=0$, which proves the claim.
\end{IEEEproof}

If the weights $w_i$ are squares $w_i=v_i^2$, which is always the case if $q=2^m$, let $g_v$ be coordinatwise multiplication by $v_{\bullet}$. Then the code $\tilde{C}= g_v(C)$ satisfy $\tilde{C} \subseteq \tilde{C}^{\perp}$ and the dual code $C'=\tilde{C}^{\perp}$ satisfy  $C'^{\perp} \subseteq C'$.

This code can be used in the Calderbank-Shor-Steane method
for constructing quantum error correcting codes.

\begin{example}(Hirzebruch Surfaces) 
The toric codes are obtained from evaluating certain rational functions in a suitable set $S$ of rational points on toric varieties, being the intersection of two ample divisors on $X$ as in \refdefn{support}.
With
\begin{align*}
F_1=&\prod_{\psi \in \F\backslash \{1\}}(e(m_1)-\psi)\ (e(m_2)-1)\\
F_2=&(e(m_1)-1) \prod_{\psi \in \F \backslash \{1\}}(e(m_2)-\psi)
\end{align*}
we have the divisors
\begin{align*}
D_1=& (F_1)_0 \sim (q-2)(V(\rho_1)+r V(\rho_4))+V(\rho_2)\\
D_2=& (F_2)_0 \sim (V(\rho_1)+r V(\rho_4))+(q-2) V(\rho_2)
\end{align*}
as
\begin{align*}
(\mathbf
e(m_{1})-\psi)_{0} \sim& (\mathbf
e(m_{1}))_{0} \sim V(\rho_1)+r V(\rho_4)\\
(e(m_{2})-\psi)_{0} \sim&  V(\rho_2)\ .
\end{align*}
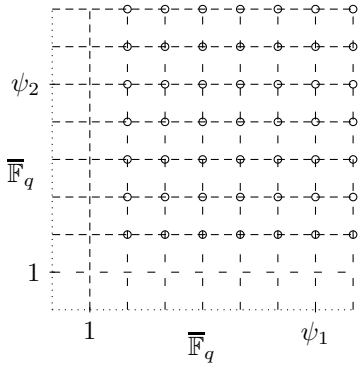
\begin{figure}
\begin{center}
\begin{tikzpicture}
\draw[dotted] (0,0) -- (0,4);
\draw[dotted] (0,0) -- (4,0);
\draw[loosely dashed](1,0)--(1,4);
\draw[loosely dashed](1.5,0)--(1.5,4);
\draw[loosely dashed](2,0)--(2,4);
\draw[loosely dashed](2.5,0)--(2.5,4);
\draw[loosely dashed](3,0)--(3,4);
\draw[loosely dashed](3.5,0)--(3.5,4);
\draw[loosely dashed](4,0)--(4,4);
\draw[loosely dashed](0,0.5)--(4,0.5);
\draw[densely dashed](0.5,0)--(0.5,4);
\draw[densely dashed](0,1)--(4,1);
\draw[densely dashed](0,1.5)--(4,1.5);
\draw[densely dashed](0,2)--(4,2);
\draw[densely dashed](0,2.5)--(4,2.5);
\draw[densely dashed](0,3)--(4,3);
\draw[densely dashed](0,3.5)--(4,3.5);
\draw[densely dashed](0,4)--(4,4);
\foreach \x in {2,3,4,5,6,7,8} \draw (1,\x/2) circle (0.05cm);
\foreach \x in {2,3,4,5,6,7,8} \draw (1.5,\x/2) circle (0.05cm);
\foreach \x in {2,3,4,5,6,7,8} \draw (2,\x/2) circle (0.05cm);
\foreach \x in {2,3,4,5,6,7,8} \draw (2.5,\x/2) circle (0.05cm);
\foreach \x in {2,3,4,5,6,7,8} \draw (3,\x/2) circle (0.05cm);
\foreach \x in {2,3,4,5,6,7,8} \draw (3.5,\x/2) circle (0.05cm);
\foreach \x in {2,3,4,5,6,7,8} \draw (4,\x/2) circle (0.05cm);
\draw (0.5,1pt) -- (0.5,-1pt) node[anchor=north] {$1$};
\draw (1pt,0.5) -- (-1pt,0.5) node[anchor=east] {$1$};
\draw (3.5,1pt) -- (3.5,-1pt) node[anchor=north] {$\psi_1$};
\draw (1pt,3) -- (-1pt,3) node[anchor=east] {$\psi_2$};
\draw (-0.1,1.8) node[anchor=east] {$\overline{\mathbb F}_q$};
\draw (2.3,-0.5) node[anchor=east] {$\overline{\mathbb F}_q$};
\end{tikzpicture}
\caption{Illustration of the supports $\vert D_1\vert$ and $\vert D_2\vert$ of the divisors $D_1$ (loosely dashed) and $D_2$ (densely dashed) and the support of the code  $S=\vert D_1\vert \cap \vert D_2 \vert$ having $(q-2)^2$ points.}\label{drawing}
\end{center}
\end{figure}
The divisors $D_1$ and $D_2$ are seen to be ample on $X$, using the intersection numbers in Table \ref{snittal} and the Nakai criterion. The support set of the code $S =\vert D_i \vert \cap \vert D_2 \vert= U_1 \cup U_2 =(\F\backslash \{1\})\times (\F\backslash \{1\}) \subseteq \F \times \F$ is realized as the intersection of the support of two ample divisors and we can apply the construction above.

\end{example}

%\nocite{*}
%\bibliographystyle{IEEEtran}
%\bibliography{IEEEabrv,IEEE-submission}

\begin{thebibliography}{10}
\providecommand{\url}[1]{#1}
\csname url@samestyle\endcsname
\providecommand{\newblock}{\relax}
\providecommand{\bibinfo}[2]{#2}
\providecommand{\BIBentrySTDinterwordspacing}{\spaceskip=0pt\relax}
\providecommand{\BIBentryALTinterwordstretchfactor}{4}
\providecommand{\BIBentryALTinterwordspacing}{\spaceskip=\fontdimen2\font plus
\BIBentryALTinterwordstretchfactor\fontdimen3\font minus
  \fontdimen4\font\relax}
\providecommand{\BIBforeignlanguage}[2]{{%
\expandafter\ifx\csname l@#1\endcsname\relax
\typeout{** WARNING: IEEEtran.bst: No hyphenation pattern has been}%
\typeout{** loaded for the language `#1'. Using the pattern for}%
\typeout{** the default language instead.}%
\else
\language=\csname l@#1\endcsname
\fi
#2}}
\providecommand{\BIBdecl}{\relax}
\BIBdecl

\bibitem{MR1749454}
J.~P. Hansen, ``Toric surfaces and error-correcting codes,'' in \emph{Coding
  theory, cryptography and related areas ({G}uanajuato, 1998)}.\hskip 1em plus
  0.5em minus 0.4em\relax Berlin: Springer, 2000, pp. 132--142.

\bibitem{MR1953195}
------, ``Toric varieties {H}irzebruch surfaces and error-correcting codes,''
  \emph{Appl. Algebra Engrg. Comm. Comput.}, vol.~13, no.~4, pp. 289--300,
  2002.

\bibitem{Calderbank19961098}
A.~Calderbank and P.~Shor, ``Good quantum error-correcting codes exist,''
  \emph{Physical Review A - Atomic, Molecular, and Optical Physics}, vol.~54,
  no.~2, pp. 1098--1105, 1996.

\bibitem{MR1450603}
P.~W. Shor, ``Fault-tolerant quantum computation,'' in \emph{37th {A}nnual
  {S}ymposium on {F}oundations of {C}omputer {S}cience ({B}urlington, {VT},
  1996)}.\hskip 1em plus 0.5em minus 0.4em\relax Los Alamitos, CA: IEEE Comput.
  Soc. Press, 1996, pp. 56--65.

\bibitem{Steane19992492}
A.~Steane, ``Enlargement of calderbank-shor-steane quantum codes,'' \emph{IEEE
  Transactions on Information Theory}, vol.~45, no.~7, pp. 2492--2495, 1999.

\bibitem{tsfasman}
A.~Ashikhmin, S.~Litsyn, and M.~Tsfasman, ``Asymptotically good quantum
  codes,'' \emph{Physical Review A - Atomic, Molecular, and Optical Physics},
  vol.~63, no.~3, pp. 1--5, 2001.

\bibitem{MR2810322}
D.~A. Cox, J.~B. Little, and H.~K. Schenck, \emph{Toric varieties}, ser.
  Graduate Studies in Mathematics.\hskip 1em plus 0.5em minus 0.4em\relax
  Providence, RI: American Mathematical Society, 2011, vol. 124.

\bibitem{MR1234037}
W.~Fulton, \emph{Introduction to toric varieties}, ser. Annals of Mathematics
  Studies.\hskip 1em plus 0.5em minus 0.4em\relax Princeton, NJ: Princeton
  University Press, 1993, vol. 131, the William H. Roever Lectures in Geometry.

\bibitem{MR922894}
T.~Oda, \emph{Convex bodies and algebraic geometry}, ser. Ergebnisse der
  Mathematik und ihrer Grenzgebiete (3) [Results in Mathematics and Related
  Areas (3)].\hskip 1em plus 0.5em minus 0.4em\relax Berlin: Springer-Verlag,
  1988, vol.~15, an introduction to the theory of toric varieties, Translated
  from the Japanese.

\bibitem{MR1866342}
S.~H. Hansen, ``Error-correcting codes from higher-dimensional varieties,''
  \emph{Finite Fields Appl.}, vol.~7, no.~4, pp. 531--552, 2001.

\bibitem{MR2142431}
D.~Joyner, ``Toric codes over finite fields,'' \emph{Appl. Algebra Engrg. Comm.
  Comput.}, vol.~15, no.~1, pp. 63--79, 2004.

\bibitem{MR2747728}
R.~Joshua and R.~Akhtar, ``Toric residue codes: {I},'' \emph{Finite Fields
  Appl.}, vol.~17, no.~1, pp. 15--50, 2011.

\bibitem{MR0435074}
P.~A. Griffiths, ``Variations on a theorem of {A}bel,'' \emph{Invent. Math.},
  vol.~35, pp. 321--390, 1976.

\bibitem{MR2464546}
E.~Kunz, \emph{Residues and duality for projective algebraic varieties}, ser.
  University Lecture Series.\hskip 1em plus 0.5em minus 0.4em\relax Providence,
  RI: American Mathematical Society, 2008, vol.~47, with the assistance of and
  contributions by David A. Cox and Alicia Dickenstein.

\bibitem{MR759943}
J.~Lipman, ``Dualizing sheaves, differentials and residues on algebraic
  varieties,'' \emph{Ast\'erisque}, no. 117, pp. ii+138, 1984.

\bibitem{MR1458757}
E.~Cattani, D.~Cox, and A.~Dickenstein, ``Residues in toric varieties,''
  \emph{Compositio Math.}, vol. 108, no.~1, pp. 35--76, 1997.

\bibitem{MR1644323}
W.~Fulton, \emph{Intersection theory}, 2nd~ed., ser. Ergebnisse der Mathematik
  und ihrer Grenzgebiete. 3. Folge. A Series of Modern Surveys in Mathematics
  [Results in Mathematics and Related Areas. 3rd Series. A Series of Modern
  Surveys in Mathematics].\hskip 1em plus 0.5em minus 0.4em\relax Berlin:
  Springer-Verlag, 1998, vol.~2.

\end{thebibliography}

\end{document}